\documentclass[12pt]{article}%
\usepackage[fleqn]{amsmath}%
\usepackage{amsfonts}%
\usepackage{amssymb}%
\usepackage{graphicx}
\providecommand{\U}[1]{\protect\rule{.1in}{.1in}}
\newtheorem{theorem}{Theorem}[section]

\newtheorem{definition}[theorem]{Definition}

\newtheorem{problem}[theorem]{Problem}

\numberwithin{equation}{section}
\newenvironment{Proof}[1][Proof]{\noindent\textbf{#1.} }{\ \rule{0.5em}{0.5em}}
\begin{document}

\title{Recovery of a potential on a quantum star graph from Weyl's matrix }
\author{Sergei A. Avdonin{\small $^{\text{1}}$}, Kira V. Khmelnytskaya$^{\text{2}}$,
Vladislav V. Kravchenko{\small $^{\text{3}}$}\\{\small $^{\text{1}}$ Department of Mathematics and Statistics, University of
Alaska, Fairbanks, AK 99775, USA}\\{\small $^{\text{2}}$ Faculty of Engineering, Autonomous University of
Queretaro, }\\{\small Cerro de las Campanas s/n, col. Las Campanas Quer\'{e}taro, Qro. C.P.
76010 M\'{e}xico}\\{\small $^{\text{3}}$ Department of Mathematics, Cinvestav, Campus
Quer\'{e}taro, }\\{\small Libramiento Norponiente \#2000, Fracc. Real de Juriquilla,
Quer\'{e}taro, Qro., 76230 M\'{e}xico}\\{\small e-mail: s.avdonin@alaska.edu, khmel@uaq.edu.mx,
vkravchenko@math.cinvestav.edu.mx}}
\maketitle

\begin{abstract}
The problem of recovery of a potential on a quantum star graph from Weyl's
matrix given at a finite number of points is considered. A method for its
approximate solution is proposed. It consists in reducing the problem to a
two-spectra inverse Sturm-Liouville problem on each edge with its posterior
solution. The overall approach is based on Neumann series of Bessel functions
(NSBF) representations for solutions of Sturm-Liouville equations, and, in
fact, the solution of the inverse problem on the quantum graph reduces to
dealing with the NSBF coefficients.

The NSBF representations admit estimates for the series remainders which are
independent of the real part of the square root of the spectral parameter.
This feature makes them especially useful for solving direct and inverse
problems requiring calculation of solutions on large intervals in the spectral
parameter. Moreover, the first coefficient of the NSBF representation alone is
sufficient for the recovery of the potential.

The knowledge of the Weyl matrix at a set of points allows one to calculate a
number of the NSBF coefficients at the end point of each edge, which leads to
approximation of characteristic functions of two Sturm-Liouville problems and
allows one to compute the Dirichlet-Dirichlet and Neumann-Dirichlet spectra on
each edge. In turn, for solving this two-spectra inverse Sturm-Liouville
problem a system of linear algebraic equations is derived for computing the
first NSBF coefficient and hence for recovering the potential. The proposed
method leads to an efficient numerical algorithm that is illustrated by a
number of numerical tests.

\end{abstract}

\section{Introduction}

The Weyl matrix of a quantum graph is in fact its Dirichlet-to-Neumann map. It
plays essential role in all aspects of study of quantum graphs, including
spectral theory \cite{Yurko2005} and controllability \cite{AvdoninKurasov2008}%
. The problem addressed in the present work has a clear physical meaning. From
a measured response at a number of frequencies of a system modelled by a
quantum graph, recover the differential operator on each edge of the graph.

A number of surveys and collections of papers on quantum graphs appeared last
years, and the first books on this topic by Berkolaiko and Kuchment
\cite{BerkolaikoKuchment}, Mugnolo \cite{Mugnolo} and Kurasov \cite{Kurasov}
contain excellent lists of references. Inverse spectral theory of network-like
structures is an important part of the rapidly developing area of applied
mathematics --- analysis on graphs. The known results in this direction
concern almost exclusively trees, i.e. graphs without cycles, see e.g.
\cite{Kurasov et al 2005, Yurko2005, Belishev Vakulenko 2006,
AvdoninKurasov2008, AvdoninLeugeringMikhaylov2010, AvdoninBell2015}.

To date, there are few papers containing numerical results for inverse
problems on graphs, all of them concern only very simple trees \cite{Belishev
Vakulenko 2006, AvdoninBelinskiyMatthews}. It is known that the problems of
space discretization of differential equations on metric graphs turn out to be
very difficult, and even the forward boundary value problems on graphs contain
a lot of numerical challenges (see, e.g. \cite{ArioliBenzi2018}).

In \cite{AvdoninKravchenko2022} we developed a new approach for solving
inverse spectral problems on compact quantum star graphs which lends itself to
efficient numerical algorithms. The approach is based on Neumann series of
Bessel functions (NSBF) representations for solutions of Sturm-Liouville
equations obtained in \cite{KNT} (see also \cite{KrBook2020}). Here we extend
this approach to the inverse problem of recovery the Schr\"{o}dinger equation
on a star graph from the values of the corresponding Weyl matrix at a finite
set of points.

The NSBF representations (see Section \ref{Sect NSBF} below) for solutions of
Sturm-Liouville equations
\[
-y^{\prime\prime}+q(x)y=\rho^{2}y,
\]
possess two remarkable features which make them especially convenient for
solving inverse problems. The remainders of the series admit $\rho
$-independent bounds for $\rho\in\mathbb{R}$, and the potential $q(x)$ can be
recovered from the very first coefficient of the series. The first feature
allows us to work with the approximate solutions on very large intervals in
$\rho$, and the second implies that computationally satisfactory results
require considering a reduced number of the terms of the series, which
eventually results in a reduced number of linear algebraic equations which
should be solved in each step.

In fact the overall developed approach reduces the solution of the inverse
problem on a graph to operations with the NSBF coefficients. In the first
step, using the given data, we compute the NSBF coefficients for solutions
satisfying the homogeneous Dirichlet and Neumann conditions at the boundary
vertices at the end point of each edge, which is associated with the common
vertex of the star graph. This first step allows us to split the problem on
the graph into separate problems on each edge. Second, the set of the
coefficients for the series representation of the solution satisfying the
homogeneous Dirichlet condition is used for computing the Dirichlet-Dirichlet
eigenvalues of the potential on each edge, while the coefficients for the
series representation of the solution satisfying the homogeneous Neumann
condition are used for computing the Neumann-Dirichlet eigenvalues. Moreover,
the first feature of the NSBF representations implies that, if necessary,
hundreds of the eigenvalues can be computed with uniform accuracy, and for
this, few coefficients of the series representations are sufficient.

Thus, the inverse problem on a graph is reduced to a two spectra inverse
Sturm-Liouville problem on each edge. Results on the uniqueness and
solvability of the two spectra problem are well known and can be found, e.g.,
in \cite{Chadan et al 1997}, \cite{LevitanInverse}, \cite{SavchukShkalikov},
\cite{Yurko2007}. For this problem we propose a method which again involves
the NSBF representations. It allows us to compute multiplier constants
\cite{Brown et al 2003} relating the Neumann-Dirichlet eigenfunctions
associated with the same eigenvalues but normalized at the opposite endpoints
of the interval. This leads to a system of linear algebraic equations for the
coefficients of their NSBF representations already for interior points of the
interval. Solving the system we find the very first coefficient, from which
the potential is recovered. It is worth mentioning that the NSBF
representations were first used for solving inverse Sturm-Liouville problems
on a finite interval in \cite{Kr2019JIIP}. Later on, the approach from
\cite{Kr2019JIIP} was improved in \cite{KrBook2020} and \cite{KT2021 IP1},
\cite{KT2021 IP2}. In these papers the system of linear algebraic equations
was obtained with the aid of the Gelfand-Levitan integral equation. In
\cite{Kr2022Completion} another approach, based on the consideration of the
eigenfunctions normalized at the opposite endpoints, was developed, and this
idea was used in \cite{AvdoninKravchenko2022}. In the present work we adjust
the same idea to the solution of the two-spectra inverse problem arising in
the final step.

In Section \ref{Sect ProblemSetting} we recall the definition of the Weyl
matrix and formulate the inverse problem. In Section \ref{Sect FundSys} we
write the system of equations in terms of the fundamental systems of solutions
on each edge, which is obtained directly from the knowledge of the Weyl
matrix. In Section \ref{Sect NSBF} we recall the NSBF representations for
solutions of the Sturm-Liouville equation and some of their relevant features.
In Section \ref{Sect SolutionDirect} we explain the solution of the direct
problem, i.e., the construction of the Weyl matrix when the potential on the
graph is known. In Section \ref{Sect Inverse} we give a detailed description
of the proposed method for the solution of the inverse problem. In Section
\ref{Sect Numerical} we discuss the numerical implementation of the method.
Finally, Section \ref{Sect Concl} contains some concluding remarks.

\section{Problem setting\label{Sect ProblemSetting}}

Let $\Omega$ denote a compact star graph consisting of $M$ edges $e_{1}%
$,...,$e_{M}$ connected at the vertex $v$. Every edge $e_{j}$ is identified
with an interval $(0,L_{j})$ of the real line in such a way that zero
corresponds to the boundary vertex $\gamma_{j}$, and the endpoint $L_{j}$
corresponds to the vertex $v$. By $\Gamma$ we denote the set of the boundary
vertices of $\Omega$, $\Gamma=\left\{  \gamma_{1},\ldots,\gamma_{M}\right\}
$. Let $q\in L_{2}(\Omega)$ be real valued, and $\lambda$ a complex number. A
continuous function $u$ defined on the graph $\Omega$ is an $M$-tuple of
functions $u_{j}\in C\left[  0,L_{j}\right]  $ satisfying the continuity
equalities at the vertex $v$: $u_{i}(v)=u_{j}(v)$ for all $1\leq i,j\leq M$.
Then $u\in C(\Omega)$.

We say that a function $u$ is a solution of the equation
\begin{equation}
-u^{\prime\prime}(x)+q(x)u(x)=\lambda u(x) \label{Schr}%
\end{equation}
on the graph $\Omega$ if besides (\ref{Schr}), the following conditions are
satisfied
\begin{equation}
u\in C(\Omega) \label{wcontinuity}%
\end{equation}
and
\begin{equation}
\sum_{j=1}^{M}\partial u_{j}(v)=0, \label{KN}%
\end{equation}
where $\partial u_{j}(v)$ denotes the derivative of $u$ at the vertex $v$
taken along the edge $e_{j}$ in the direction outward the vertex. The sum in
the equality (\ref{KN}), which is known as the Kirchhoff-Neumann condition, is
taken over all the edges $e_{j}$, $j=1,\ldots,M$.

Let $\lambda\notin\mathbb{R}$ and $w_{i}$ be the so-called Weyl solution, that
is, a solution of (\ref{Schr}) on $\Omega$, satisfying the initial conditions
at the boundary vertices
\[
w_{i}(\gamma_{i})=1\quad\text{and}\quad w_{i}(\gamma_{j})=0\text{ for all
}j\neq i.
\]

\begin{definition}
The $M\times M$ matrix-function $\mathbf{M}(\lambda)$, $\lambda\notin
\mathbb{R}$, consisting of the elements $\mathbf{M}_{ij}(\lambda)=\partial
w_{i}(\gamma_{j})$, $i,j=1,\ldots,M$ is called the Weyl matrix.
\end{definition}

In fact, for a fixed value of $\lambda$, the Weyl matrix is the
Dirichlet-to-Neumann map on the quantum graph because, if $u$ is a solution of
(\ref{Schr}) satisfying the Dirichlet condition at the boundary vertices
$u(\gamma,\lambda)=f(\lambda)$, then $\partial u(\gamma,\lambda)=\mathbf{M}%
(\lambda)f(\lambda)$, $\lambda\notin\mathbb{R}$.

The problem we consider in the present paper can be formulated as follows.

\begin{problem}
Given the Weyl matrix at a finite number of points $\lambda_{k}$,
$k=1,\ldots,m$, approximate the potential $q(x)$.
\end{problem}

When the Weyl matrix or even its main diagonal is known everywhere, the
potential $q(x)$ is determined uniquely (see, e.g., \cite{Yurko2005}). The
knowledge of the Weyl matrix at a finite number of points allows one to
recover the potential $q(x)$ only approximately.

\section{Fundamental system of solutions and the Weyl
matrix\label{Sect FundSys}}

In order to reformulate the inverse problem in terms suitable for its
numerical solution, let us introduce for each edge $e_{i}$ a corresponding
fundamental system of solutions. By $\varphi_{i}(\rho,x)$ and $S_{i}(\rho,x)$
we denote the solutions of the equation
\begin{equation}
-y^{\prime\prime}(x)+q_{i}(x)y(x)=\rho^{2}y(x),\quad x\in(0,L_{i})
\label{Schri}%
\end{equation}
satisfying the initial conditions
\[
\varphi_{i}(\rho,0)=1,\quad\varphi_{i}^{\prime}(\rho,0)=0,
\]%
\[
S_{i}(\rho,0)=0,\quad S_{i}^{\prime}(\rho,0)=1.
\]
Here $q_{i}(x)$ is the potential $q(x)$ on the edge $e_{i}$ and $\rho
=\sqrt{\lambda}$, $\operatorname{Im}\rho\geq0$. Then the Weyl solution $w_{i}$
has the form%
\[
w_{i}(\rho,x)=\varphi_{i}(\rho,x)+\mathbf{M}_{ii}(\rho^{2})S_{i}(\rho
,x)\quad\text{on the edge }e_{i}%
\]
and
\[
w_{i}(\rho,x)=\mathbf{M}_{ij}(\rho^{2})S_{j}(\rho,x)\quad\text{on every edge
}e_{j},\quad j\neq i.
\]

Thus, the knowledge of the Weyl matrix at a point $\rho_{k}^{2}$ implies the
knowledge of such constants $\mathbf{M}_{ij}(\rho_{k}^{2})$ that for all
$i=1,\ldots,M$, the equalities are valid%
\begin{equation}
\varphi_{i}(\rho_{k},L_{i})+\mathbf{M}_{ii}(\rho_{k}^{2})S_{i}(\rho_{k}%
,L_{i})=\mathbf{M}_{ij}(\rho_{k}^{2})S_{j}(\rho_{k},L_{j})\quad\text{for all
}j\neq i \label{continuity i}%
\end{equation}
and
\begin{equation}
\varphi_{i}^{\prime}(\rho_{k},L_{i})+\sum_{j=1}^{M}\mathbf{M}_{ij}(\rho
_{k}^{2})S_{j}^{\prime}(\rho_{k},L_{j})=0. \label{KN i}%
\end{equation}

The solutions $\varphi_{i}(\rho,x)$ and $S_{i}(\rho,x)$ admit very convenient
series representations, which we discuss in the next section.

\section{Neumann series of Bessel functions representations\label{Sect NSBF}}

In analysis, according to \cite[Chapter XVI]{Watson}, \textquotedblleft Any
series of the type
\[
\sum_{n=0}^{\infty}a_{n}J_{\nu+n}(z)
\]
is called a Neumann series, although in fact Neumann considered only
the\ special type of series for which $\nu$ is an integer; the investigation
of the more general series is due to Gegenbauer\textquotedblright. The papers
by C. G. Neumann and L. B. Gegenbauer date from 1867 and 1877, respectively.
Since then the Neumann series of Bessel functions (NSBF) were studied in
numerous publications (see \cite{Watson}, \cite{Wilkins} and the recent
monograph on the subject \cite{Baricz et al Book} and references therein). In
\cite{KNT} it was shown that the solutions $\varphi_{i}(\rho,x)$, $S_{i}%
(\rho,x)$ and their derivatives with respect to $x$ admit the NSBF
representations which possess certain unique features which make them
especially useful for solving direct and inverse spectral problems.

\begin{theorem}
[\cite{KNT}]\label{Th NSBF} The solutions $\varphi_{i}(\rho,x)$ and
$S_{i}(\rho,x)$ of (\ref{Schri}) and their derivatives with respect to $x$
admit the following series representations
\begin{align}
\varphi_{i}(\rho,x)  &  =\cos\left(  \rho x\right)  +\sum_{n=0}^{\infty
}(-1)^{n}g_{i,n}(x)\mathbf{j}_{2n}(\rho x),\label{phiNSBF}\\
S_{i}(\rho,x)  &  =\frac{\sin\left(  \rho x\right)  }{\rho}+\frac{1}{\rho}%
\sum_{n=0}^{\infty}(-1)^{n}s_{i,n}(x)\mathbf{j}_{2n+1}(\rho x),\label{S}\\
\varphi_{i}^{\prime}(\rho,x)  &  =-\rho\sin\left(  \rho x\right)  +\frac
{\cos\left(  \rho x\right)  }{2}\int_{0}^{x}q_{i}(t)\,dt+\sum_{n=0}^{\infty
}(-1)^{n}\gamma_{i,n}(x)\mathbf{j}_{2n}(\rho x),\label{phiprimeNSBF}\\
S_{i}^{\prime}(\rho,x)  &  =\cos\left(  \rho x\right)  +\frac{\sin\left(  \rho
x\right)  }{2\rho}\int_{0}^{x}q_{i}(t)\,dt+\frac{1}{\rho}\sum_{n=0}^{\infty
}(-1)^{n}\sigma_{i,n}(x)\mathbf{j}_{2n+1}(\rho x), \label{Sprime}%
\end{align}
where $\mathbf{j}_{k}(z)$ stands for the spherical Bessel function of order
$k$, $\mathbf{j}_{k}(z):=\sqrt{\frac{\pi}{2z}}J_{k+\frac{1}{2}}(z)$ (see,
e.g., \cite{AbramowitzStegunSpF}). The coefficients $g_{i,n}(x)$, $s_{i,n}%
(x)$, $\gamma_{i,n}(x)$ and $\sigma_{i,n}(x)$ can be calculated following a
simple recurrent integration procedure (see \cite{KNT} or \cite[Sect.
9.4]{KrBook2020}), starting with
\begin{align}
g_{i,0}(x)  &  =\varphi_{i}(0,x)-1,\quad s_{i,0}(x)=3\left(  \frac{S_{i}%
(0,x)}{x}-1\right)  ,\label{beta0}\\
\gamma_{i,0}(x)  &  =g_{i,0}^{\prime}(x)-\frac{1}{2}\int_{0}^{x}%
q_{i}(t)\,dt,\quad\sigma_{i,0}(x)=\frac{s_{i,0}(x)}{x}+s_{i,0}^{\prime
}(x)-\frac{3}{2}\int_{0}^{x}q_{i}(t)\,dt.\nonumber
\end{align}
For every $\rho\in\mathbb{C}$ all the series converge pointwise. For every
$x\in\left[  0,L_{i}\right]  $ the series converge uniformly on any compact
set of the complex plane of the variable $\rho$, and the remainders of their
partial sums admit estimates independent of $\operatorname{Re}\rho$.
\end{theorem}

This last feature of the series representations (the independence of
$\operatorname{Re}\rho$ of the estimates for the remainders) is of crucial
importance for what follows. In particular, it means that for $S_{i,N}%
(\rho,x):=\frac{\sin\left(  \rho x\right)  }{\rho}+\frac{1}{\rho}\sum
_{n=0}^{N}(-1)^{n}s_{i,n}(x)\mathbf{j}_{2n+1}(\rho x)$ (and analogously for
$\varphi_{i,N}(\rho,x)$) the estimate holds
\begin{equation}
\left\vert S_{i}(\rho,x)-S_{i,N}(\rho,x)\right\vert <\varepsilon_{i,N}(x)
\label{estim S}%
\end{equation}
for all $\rho\in\mathbb{R}$, where $\varepsilon_{i,N}(x)$ is a positive
function tending to zero as $N\rightarrow\infty$. Roughly speaking, the
approximate solution $S_{i,N}(\rho,x)$ approximates the exact one equally well
for small and for large values of $\rho\in\mathbb{R}$. This is especially
convenient when considering direct and inverse spectral problems that requires
operating on a large range of the parameter $\rho$. This unique feature of the
series representations (\ref{phiNSBF})-(\ref{Sprime}) is due to the fact that
they originate from an exact Fourier-Legendre series representation of the
integral kernel of the transmutation operator \cite{KNT}, \cite[Sect.
9.4]{KrBook2020} (for the theory of transmutation operators we refer to
\cite{LevitanInverse}, \cite{Marchenko}, \cite{SitnikShishkina Elsevier},
\cite{Yurko2007}.

Moreover, for a fixed $z$ the numbers $\mathbf{j}_{k}(z)$ rapidly decrease as
$k\rightarrow\infty$, see, e.g., \cite[(9.1.62)]{AbramowitzStegunSpF}. Hence,
the convergence rate of the series for any fixed $\rho$ is, in fact,
exponential. More detailed estimates for the series remainders depending on
the regularity of the potential can be found in \cite{KNT}.

Note that formulas (\ref{beta0}) indicate that the potential $q_{i}(x)$ can be
recovered from the first coefficients of the series (\ref{phiNSBF}) or
(\ref{S}). Indeed, we have
\begin{equation}
q_{i}(x)=\frac{g_{i,0}^{\prime\prime}(x)}{g_{i,0}(x)+1} \label{qi from g0}%
\end{equation}
and
\begin{equation}
q_{i}(x)=\frac{\left(  xs_{i,0}(x)\right)  ^{\prime\prime}}{xs_{i,0}(x)+3x}.
\label{qi from s0}%
\end{equation}

\section{Solution of the direct problem\label{Sect SolutionDirect}}

Here we briefly explain how in terms of the solutions $\varphi_{i}(\rho,x)$
and $S_{i}(\rho,x)$ the construction of the Weyl matrix can be performed.
Given $q(x)$, fix $\rho$. To find the first row of the matrix $\mathbf{M}%
(\rho^{2})$, that is the entries $\mathbf{M}_{1j}(\rho^{2})$, we consider the
corresponding equalities at the common vertex $v$:
\[
\varphi_{1}(\rho,L_{1})+\mathbf{M}_{11}(\rho^{2})S_{1}(\rho,L_{1}%
)=\mathbf{M}_{12}(\rho^{2})S_{2}(\rho,L_{2}),
\]%
\[
\mathbf{M}_{12}(\rho^{2})S_{2}(\rho,L_{2})=\mathbf{M}_{13}(\rho^{2})S_{3}%
(\rho,L_{3}),
\]%
\[
\ldots
\]%
\[
\mathbf{M}_{1M-1}(\rho^{2})S_{M-1}(\rho,L_{M-1})=\mathbf{M}_{1M}(\rho
^{2})S_{M}(\rho,L_{M}),
\]%
\[
\varphi_{1}^{\prime}(\rho,L_{1})+\mathbf{M}_{11}(\rho^{2})S_{1}^{\prime}%
(\rho,L_{1})+\mathbf{M}_{12}(\rho^{2})S_{2}^{\prime}(\rho,L_{2})+\ldots
+\mathbf{M}_{1M}(\rho^{2})S_{M}^{\prime}(\rho,L_{M})=0,
\]
which can be written in the form of a system of linear algebraic equations%
\[
\mathbf{A}(\rho)\overrightarrow{\mathbf{M}}_{1}(\rho^{2})=
\]%
\begin{equation}
\left(
\begin{array}
[c]{cccccc}%
S_{1}(L_{1}) & -S_{2}(L_{2}) & 0 & 0 & \ldots & 0\\
0 & S_{2}(L_{2}) & -S_{3}(L_{3}) & 0 & \ldots & 0\\
&  & \ldots &  &  & \\
0 & \ldots &  & 0 & S_{M-1}(L_{M-1}) & -S_{M}(L_{M})\\
S_{1}^{\prime}(L_{1}) & S_{2}^{\prime}(L_{2}) & \ldots &  & S_{M-1}^{\prime
}(L_{M-1}) & S_{M}^{\prime}(L_{M})
\end{array}
\right)  \left(
\begin{array}
[c]{c}%
\mathbf{M}_{11}\\
\mathbf{M}_{12}\\
\vdots\\
\mathbf{M}_{1M-1}\\
\mathbf{M}_{1M}%
\end{array}
\right)  =\left(
\begin{array}
[c]{c}%
-\varphi_{1}(L_{1})\\
0\\
\vdots\\
0\\
-\varphi_{1}^{\prime}(L_{1})
\end{array}
\right)  , \label{AM1}%
\end{equation}
where for the sake of space we omitted the dependence on $\rho$.

It is easy to verify that all subsequent rows of the Weyl matrix can be
computed by solving a system of linear algebraic equations with the same
matrix
\[
\mathbf{A}(\rho)=\left(
\begin{array}
[c]{cccccc}%
S_{1}(\rho,L_{1}) & -S_{2}(\rho,L_{2}) & 0 & 0 & \ldots & 0\\
0 & S_{2}(\rho,L_{2}) & -S_{3}(\rho,L_{3}) & 0 & \ldots & 0\\
&  & \ldots &  &  & \\
0 & \ldots &  & 0 & S_{M-1}(\rho,L_{M-1}) & -S_{M}(\rho,L_{M})\\
S_{1}^{\prime}(\rho,L_{1}) & S_{2}^{\prime}(\rho,L_{2}) & \ldots &  &
S_{M-1}^{\prime}(\rho,L_{M-1}) & S_{M}^{\prime}(\rho,L_{M})
\end{array}
\right)  ,
\]
but with a different right-hand side. Namely,
\begin{equation}
\mathbf{A}(\rho)\overrightarrow{\mathbf{M}}_{i}(\rho^{2})=\left(
\begin{array}
[c]{c}%
0\\
\vdots\\
0\\
\varphi_{i}(\rho,L_{i})\\
-\varphi_{i}(\rho,L_{i})\\
0\\
\vdots\\
-\varphi_{i}^{\prime}(\rho,L_{i})
\end{array}
\right)  \quad\text{for }i=2,\ldots,M-1, \label{AMi}%
\end{equation}
where the first $\varphi_{i}(\rho,L_{i})$ appears at $\left(  i-1\right)  $-th
position, and for the last row we have
\begin{equation}
\mathbf{A}(\rho)\overrightarrow{\mathbf{M}}_{M}(\rho^{2})=\left(
\begin{array}
[c]{c}%
0\\
\vdots\\
0\\
\varphi_{M}(\rho,L_{M})\\
-\varphi_{M}^{\prime}(\rho,L_{M})
\end{array}
\right)  . \label{AMM}%
\end{equation}

\section{Solution of the inverse problem\label{Sect Inverse}}

Here we assume the Weyl matrix to be known at a finite number of points
$\rho_{k}^{2}\notin\mathbb{R}$, $k=1,\ldots,m$. From this information we look
to recover the potential $q_{i}(x)$ on each edge of the star graph. The
proposed method for solving this problem consists of several steps.

First, we compute a number of the constants $\left\{  g_{i,n}(L_{i})\right\}
_{n=0}^{N}$ and $\left\{  s_{i,n}(L_{i})\right\}  _{n=0}^{N}$ for every
$i=1,\ldots,M$, that is, the values of the coefficients from (\ref{phiNSBF})
and (\ref{S}) at the endpoint of the edge $e_{i}$. This first step allows us
to split the problem on the graph and reduce it to the problems on the edges.
Indeed, the knowledge of the coefficients $\left\{  g_{i,n}(L_{i})\right\}
_{n=0}^{N}$ and $\left\{  s_{i,n}(L_{i})\right\}  _{n=0}^{N}$ allows us in the
second step to compute the Dirichlet-Dirichlet and Neumann-Dirichlet spectra
for the potential $q_{i}(x)$, $x\in\left[  0,L_{i}\right]  $, thus obtaining a
two-spectra inverse problem for each component of the potential $q(x)$.

In the third step the two-spectra inverse problem is solved with the aid of
the representation (\ref{phiNSBF}) and an analogous series representation for
the solution $T_{i}(\rho,x)$ of (\ref{Schri}) satisfying the initial
conditions at the endpoint $L_{i}$:
\begin{equation}
T_{i}(\rho,L_{i})=0,\quad T_{i}^{\prime}(\rho,L_{i})=1. \label{psi init}%
\end{equation}
The two-spectra inverse problem is reduced to a system of linear algebraic
equations, from which we obtain the coefficient $g_{i,0}(x)$. Finally, the
potential $q_{i}(x)$ is calculated from (\ref{qi from g0}).

\subsection{Calculation of coefficients $\left\{  g_{i,n}(L_{i})\right\}
_{n=0}^{N}$ and $\left\{  s_{i,n}(L_{i})\right\}  _{n=0}^{N}$}

The knowledge of the Weyl matrix at a point $\rho_{k}^{2}$ means that we have
equalities (\ref{AM1})-(\ref{AMM}) valid for $\rho=\rho_{k}$. In fact we will
use only the equations which arise from the continuity condition
(\ref{wcontinuity}) and not those which arise from the Kirchhoff-Neumann
condition (\ref{KN}). Hence we work with the solutions $\varphi_{i}(\rho,x)$
and $S_{i}(\rho,x)$ and not with their derivatives. Thus, for every $\rho_{k}$
from the set, using the representations (\ref{phiNSBF}) and (\ref{S}) we have
the equations
\begin{align}
&  \rho_{k}\sum_{n=0}^{\infty}(-1)^{n}g_{i,n}(L_{i})\mathbf{j}_{2n}(\rho
_{k}L_{i})+\mathbf{M}_{ii}(\rho_{k}^{2})\sum_{n=0}^{\infty}(-1)^{n}%
s_{i,n}(L_{i})\mathbf{j}_{2n+1}(\rho_{k}L_{i})\nonumber\\
&  -\mathbf{M}_{ii+1}(\rho_{k}^{2})\sum_{n=0}^{\infty}(-1)^{n}s_{i+1,n}%
(L_{i+1})\mathbf{j}_{2n+1}(\rho_{k}L_{i+1})\nonumber\\
&  =\mathbf{M}_{ii+1}(\rho_{k}^{2})\sin(\rho_{k}L_{i+1})-\rho_{k}\cos(\rho
_{k}L_{i})-\mathbf{M}_{ii}(\rho_{k}^{2})\sin(\rho_{k}L_{i}),\quad\text{for
}i=1,\ldots,M, \label{type1}%
\end{align}
where for $i=M$ we replace $i+1$ by $1$ (the cyclic rule), as well as the
equations%
\begin{align}
&  \mathbf{M}_{ij}(\rho_{k}^{2})\sum_{n=0}^{\infty}(-1)^{n}s_{j,n}%
(L_{j})\mathbf{j}_{2n+1}(\rho_{k}L_{j})-\mathbf{M}_{ij+1}(\rho_{k}^{2}%
)\sum_{n=0}^{\infty}(-1)^{n}s_{j+1,n}(L_{j+1})\mathbf{j}_{2n+1}(\rho
_{k}L_{j+1})\nonumber\\
&  =\mathbf{M}_{ij+1}(\rho_{k}^{2})\sin(\rho_{k}L_{j+1})-\mathbf{M}_{ij}%
(\rho_{k}^{2})\sin(\rho_{k}L_{j}), \label{type2}%
\end{align}
where $j=1,\ldots,M$, $j\neq i$, $j+1\neq i$, and again, for $j=M$ we replace
$j+1$ by $1$.

To compute the finite sets of the coefficients $\left\{  g_{i,n}%
(L_{i})\right\}  _{n=0}^{N}$ and $\left\{  s_{i,n}(L_{i})\right\}  _{n=0}^{N}$
for all $i=1,\ldots,M$ we chose the following strategy. Fix $i$, and for each
$\rho_{k}^{2}$ consider the continuity conditions corresponding to the Weyl
solution $w_{i}$ or, in other words, to the $i$-th row of the Weyl matrix. We
have then one equation of the form (\ref{type1}) and $M-2$ equations of the
form (\ref{type2}). Considering these $M-1$ equations for every $\rho_{k}^{2}%
$, we obtain $m(M-1)$ equations for $(M+1)(N+1)$ unknowns. Here the unknowns
are $M$ sets of the coefficients $\left\{  s_{j,n}(L_{i})\right\}  _{n=0}^{N}$
(for $j=1,\ldots,M$) and one set of the coefficients $\left\{  g_{i,n}%
(L_{i})\right\}  _{n=0}^{N}$. Thus, we need at least $m=\left\lceil
\frac{(M+1)(N+1)}{M-1}\right\rceil $ points $\rho_{k}^{2}$ at which the Weyl
matrix is known. Here $\left\lceil x\right\rceil $ denotes the least integer
greater than or equal to $x$. In particular, the obtained system of equations
allows us to compute $\left\{  g_{i,n}(L_{i})\right\}  _{n=0}^{N}$ and
$\left\{  s_{i,n}(L_{i})\right\}  _{n=0}^{N}$ for the fixed $i$. Thus in
total, we consider $M$ linear algebtraic systems of this kind to compute all
sets of the coefficients $\left\{  g_{i,n}(L_{i})\right\}  _{n=0}^{N}$ and
$\left\{  s_{i,n}(L_{i})\right\}  _{n=0}^{N}$ for all $i=1,\ldots,M$.

An important observation however consists in the fact that in order to find
the coefficients $\left\{  g_{i,n}(L_{i})\right\}  _{n=0}^{N}$ and $\left\{
s_{i,n}(L_{i})\right\}  _{n=0}^{N}$, $i=1,\ldots,M$, there is no need to
consider the equations of the form (\ref{type2}) or at least all such
equations. In fact, it is enough to consider equation (\ref{type1}) alone,
which means that to recover the potential on the edge $e_{i}$, we need to know
the main diagonal entry $\mathbf{M}_{ii}(\rho_{k}^{2})$ as well as some
$\mathbf{M}_{ij}(\rho_{k}^{2})$ for one $j\neq i$, which can be $\mathbf{M}%
_{ii+1}(\rho_{k}^{2})$ (for $i=M$, $i+1$ is replaced by $1$). In this case,
for each edge $e_{i}$, from equations of the form (\ref{type1}), which we have
for each $\rho_{k}$, we compute $\left\{  g_{i,n}(L_{i})\right\}  _{n=0}^{N}$,
$\left\{  s_{i,n}(L_{i})\right\}  _{n=0}^{N}$ and $\left\{  s_{i+1,n}%
(L_{i+1})\right\}  _{n=0}^{N}$, that is, $3(N+1)$ unknowns. For this, the
knowledge of the Weyl matrix (or of its main diagonal plus one more entry from
each row) is required at least at $3(N+1)$ points $\rho_{k}$.

We show below that eventually the accuracy of the recovered potential is
comparable when one uses all equations of the form (\ref{type2}), part of them
or none. The number of equations of the form (\ref{type2}) used in
computations we will denote by $M_{k}$. Thus, if no such equation is used
$M_{k}=0$ (that is, two elements from each row of the Weyl matrix:
$\mathbf{M}_{ii}(\rho_{k}^{2})$ and $\mathbf{M}_{ii+1}(\rho_{k}^{2})$ are
used), while $M_{k}=M-2$ means that all equations of the form (\ref{type2})
(and all elements of the Weyl matrix) are used.

\subsection{Reduction to the two-spectra inverse problem on the edge
\label{Subsect reduction to two spectra}}

The first step, as described in the previous subsection, reduces the problem
on the graph to $M$ separate problems on the edges. Thus, consider an edge
$e_{i}$ for which at this stage we have computed the coefficients $\left\{
g_{i,n}(L_{i})\right\}  _{n=0}^{N}$ and $\left\{  s_{i,n}(L_{i})\right\}
_{n=0}^{N}$. Now we use them to compute the Dirichlet-Dirichlet and
Neumann-Dirichlet spectra for the potential $q_{i}(x)$, $x\in\left[
0,L_{i}\right]  $. This is done with the aid of the approximate solutions
evaluated at the end point%
\[
\varphi_{i,N}(\rho,L_{i})=\cos\left(  \rho L_{i}\right)  +\sum_{n=0}%
^{N}(-1)^{n}g_{i,n}(L_{i})\mathbf{j}_{2n}(\rho L_{i})
\]
and
\[
S_{i,N}(\rho,L_{i})=\frac{\sin\left(  \rho L_{i}\right)  }{\rho}+\frac{1}%
{\rho}\sum_{n=0}^{N}(-1)^{n}s_{i,n}(L_{i})\mathbf{j}_{2n+1}(\rho L_{i}).
\]
Indeed, since $S_{i}(\rho,x)$ (as well as $S_{i,N}(\rho,x)$) satisfies the
Dirichlet condition at the origin, zeros of $S_{i}(\rho,L_{i})$ are precisely
square roots of the Dirichlet-Dirichlet eigenvalues. That is, $S_{i}%
(\rho,L_{i})$ is the characteristic function of the Sturm-Liouville problem
\begin{equation}
-y^{\prime\prime}+q_{i}(x)y=\lambda y,\quad x\in(0,L_{i}), \label{SLi}%
\end{equation}%
\begin{equation}
y(0)=y(L_{i})=0, \label{DD cond}%
\end{equation}
and its zeros coincide with the numbers $\left\{  \mu_{i,k}\right\}
_{k=1}^{\infty}$, such that $\mu_{i,k}^{2}$ are the eigenvalues of the problem
(\ref{SLi}), (\ref{DD cond}). In turn, zeros of the function $S_{i,N}%
(\rho,L_{i})$ approximate zeros of $S_{i}(\rho,L_{i})$ (see Theorem
\ref{Th closeness of zeros} below). Thus, the singular numbers $\left\{
\mu_{i,k}\right\}  _{k=1}^{\infty}$ are approximated by zeros of the function
$S_{i,N}(\rho,L_{i})$.

The same reasoning is valid for the function $\varphi_{i,N}(\rho,L_{i})$,
whose zeros approximate the singular numbers $\left\{  \nu_{i,k}\right\}
_{k=1}^{\infty}$, which are the square roots of the eigenvalues of the
Sturm-Liouville problem for (\ref{SLi}) subject to the boundary conditions
\begin{equation}
y^{\prime}(0)=y(L_{i})=0. \label{ND cond}%
\end{equation}
\qquad

The following statement is valid.

\begin{theorem}
\label{Th closeness of zeros}For any $\varepsilon>0$ there exists such
$N\in\mathbb{N}$ that all zeros of the functions $S_{i}(\rho,L_{i})$ and
$\varphi_{i}(\rho,L_{i})$ are approximated by corresponding zeros of the
functions $S_{i,N}(\rho,L_{i})$ and $\varphi_{i,N}(\rho,L_{i})$, respectively,
with errors uniformly bounded by $\varepsilon$. Moreover, $S_{i,N}(\rho
,L_{i})$ and $\varphi_{i,N}(\rho,L_{i})$, have no other zeros.
\end{theorem}

\begin{Proof}
The proof of this statement is completely analogous to the proof of
Proposition 7.1 in \cite{KT2015JCAM} and consists in the use of properties of
characteristic functions of regular Sturm-Liouville problems and application
of the Rouch\'{e} theorem.
\end{Proof}

Below, in Section \ref{Sect Numerical} we show that indeed, even for
relatively small $N$, computing zeros of the functions $S_{i,N}(\rho,L_{i})$
and $\varphi_{i,N}(\rho,L_{i})$ one obtains hundreds of the
Dirichlet-Dirichlet and Neumann-Dirichlet eigenvalues computed with remarkably
uniform accuracy. Thus, on every edge $e_{i}$ we obtain the classical inverse
problem of recovering the potential $q_{i}(x)$ from two spectra, which is
considered in the next step.

\subsection{Solution of two-spectra inverse problem}

At this stage we dispose of two finite sequences of singular numbers $\left\{
\mu_{i,k}\right\}  _{k=1}^{K_{D}}$ and $\left\{  \nu_{i,k}\right\}
_{k=1}^{K_{N}}$ which are square roots of the eigenvalues of problems
(\ref{SLi}), (\ref{DD cond}) and (\ref{SLi}), (\ref{ND cond}), respectively,
as well as of two sequences of numbers $\left\{  s_{i,n}(L_{i})\right\}
_{n=0}^{N}$ and $\left\{  g_{i,n}(L_{i})\right\}  _{n=0}^{N}$, which are the
values of the coefficients from (\ref{S}) and (\ref{phiNSBF}) at the endpoint.

Let us consider the solution $T_{i}(\rho,x)$ of equation (\ref{SLi})
satisfying the initial conditions at $L_{i}$:
\[
T_{i}(\rho,L_{i})=0,\quad T_{i}^{\prime}(\rho,L_{i})=1.
\]
Analogously to the solution (\ref{S}), the solution $T_{i}(\rho,x)$ admits the
series representation%
\begin{equation}
T_{i}(\rho,x)=\frac{\sin\left(  \rho\left(  x-L_{i}\right)  \right)  }{\rho
}+\frac{1}{\rho}\sum_{n=0}^{\infty}(-1)^{n}t_{i,n}(x)\mathbf{j}_{2n+1}%
(\rho\left(  x-L_{i}\right)  ), \label{Ti}%
\end{equation}
where $t_{i,n}\left(  x\right)  $ are corresponding coefficients, analogous to
$s_{i,n}\left(  x\right)  $ from (\ref{S}).

Note that for $\rho=\nu_{i,k}$ the solutions $\varphi_{i}(\nu_{i,k},x)$ and
$T_{i}(\nu_{i,k},x)$ are linearly dependent because both are eigenfunctions of
problem (\ref{SLi}), (\ref{ND cond}). Hence there exist such real constants
$\beta_{i,k}\neq0$, that
\begin{equation}
\varphi_{i}(\nu_{i,k},x)=\beta_{i,k}T_{i}(\nu_{i,k},x). \label{phi=T}%
\end{equation}
Moreover, these multiplier constants can be easily calculated by recalling
that $\varphi_{i}(\nu_{i,k},0)=1$. Thus,%
\begin{align}
\frac{1}{\beta_{i,k}}  &  =T_{i}(\nu_{i,k},0)\approx T_{i,N}(\nu
_{i,k},0)\nonumber\\
&  =-\frac{\sin\left(  \nu_{i,k}L_{i}\right)  }{\nu_{i,k}}-\frac{1}{\nu_{i,k}%
}\sum_{n=0}^{N}(-1)^{n}t_{i,n}(0)\mathbf{j}_{2n+1}(\nu_{i,k}L_{i}),
\label{betaik}%
\end{align}
where we took into account that the spherical Bessel functions of odd order
are odd. The coefficients $\left\{  t_{i,n}(0)\right\}  _{n=0}^{N}$ are
computed with the aid of the singular numbers $\left\{  \mu_{i,k}\right\}
_{k=1}^{K_{D}}$ as follows. Since the functions $T_{i}(\mu_{i,k},x)$,
$k=1,2,\ldots$ are eigenfunctions of the problem (\ref{SLi}), (\ref{DD cond}),
we have that $T_{i}(\mu_{i,k},0)=0$ and hence%
\[
\sum_{n=0}^{\infty}(-1)^{n}t_{i,n}(0)\mathbf{j}_{2n+1}(\mu_{i,k}L_{i}%
)=-\sin\left(  \mu_{i,k}L_{i}\right)  ,\quad k=1,2,\cdots.
\]
This leads to a system of linear algebraic equations for computing the
coefficients $\left\{  t_{i,n}(0)\right\}  _{n=0}^{N}$, which has the form%
\[
\sum_{n=0}^{N}(-1)^{n}t_{i,n}(0)\mathbf{j}_{2n+1}(\mu_{i,k}L_{i})=-\sin\left(
\mu_{i,k}L_{i}\right)  ,\quad k=1,\cdots,K_{D}.
\]
Now, having computed $\left\{  t_{i,n}(0)\right\}  _{n=0}^{N}$, we compute the
multiplier constants $\left\{  \beta_{i,k}\right\}  _{k=1}^{K_{N}}$ from
(\ref{betaik}).

Next, we use equation (\ref{phi=T}) for constructing a system of linear
algebraic equations for the coefficients $g_{i,n}(x)$ and $t_{i,n}\left(
x\right)  $. Indeed, equation (\ref{phi=T}) can be written in the form%
\begin{align*}
&  \sum_{n=0}^{\infty}(-1)^{n}g_{i,n}(x)\mathbf{j}_{2n}(\nu_{i,k}%
x)-\frac{\beta_{i,k}}{\nu_{i,k}}\sum_{n=0}^{\infty}(-1)^{n}t_{i,n}%
(x)\mathbf{j}_{2n+1}(\nu_{i,k}\left(  x-L_{i}\right)  )\\
&  =\frac{\beta_{i,k}}{\nu_{i,k}}\sin\left(  \nu_{i,k}\left(  x-L_{i}\right)
\right)  -\cos\left(  \nu_{i,k}x\right)  .
\end{align*}

We have as many of such equations as many Neumann-Dirichlet singular numbers
$\nu_{i,k}$ are computed. For computational purposes we choose some natural
number $N_{c}$ - the number of the coefficients $g_{i,n}(x)$ and
$t_{i,n}\left(  x\right)  $ to be computed. More precisely, we choose a
sufficiently dense set of points $x_{m}\in(0,L_{i})$ and at every $x_{m}$
consider the equations%
\begin{align*}
&  \sum_{n=0}^{N_{c}}(-1)^{n}g_{i,n}(x_{m})\mathbf{j}_{2n}(\nu_{i,k}%
x_{m})-\frac{\beta_{i,k}}{\nu_{i,k}}\sum_{n=0}^{N_{c}}(-1)^{n}t_{i,n}%
(x_{m})\mathbf{j}_{2n+1}(\nu_{i,k}\left(  x_{m}-L_{i}\right)  )\\
&  =\frac{\beta_{i,k}}{\nu_{i,k}}\sin\left(  \nu_{i,k}\left(  x_{m}%
-L_{i}\right)  \right)  -\cos\left(  \nu_{i,k}x_{m}\right)  ,\quad
k=1,\ldots,K_{N}.
\end{align*}
Solving this system of equations we find $g_{i,0}(x_{m})$ and consequently
$g_{i,0}(x)$ at a sufficiently dense set of points of the interval $(0,L_{i}%
)$. Finally, with the aid of (\ref{qi from g0}) we compute $q_{i}(x)$.

Schematically the proposed method for solving the inverse problem on a quantum
star graph is presented in the following diagram.%
\[
\mathbf{M}(\rho_{k}^{2}),\,k=1,\ldots,m\quad\overset{(1)}{\Longrightarrow
}\quad\left\{  g_{i,n}(L_{i}),\ s_{i,n}(L_{i})\right\}  _{n=0}^{N}%
\quad\overset{(2)}{\Longrightarrow}\quad\text{two spectra }\left\{  \mu
_{i,k}\right\}  _{k=1}^{K_{D}},\ \left\{  \nu_{i,k}\right\}  _{k=1}^{K_{N}}%
\]%
\[
\left\{  \mu_{i,k}\right\}  _{k=1}^{K_{D}}\quad\overset{(3)}{\Longrightarrow
}\quad\left\{  t_{i,n}(0)\right\}  _{n=0}^{N}%
\]%
\[
\left\{  t_{i,n}(0)\right\}  _{n=0}^{N},\,\left\{  \nu_{i,k}\right\}
_{k=1}^{K_{N}}\quad\overset{(4)}{\Longrightarrow}\quad\left\{  \beta
_{i,k}\right\}  _{k=1}^{K_{N}}%
\]%
\[
\left\{  \nu_{i,k},\beta_{i,k}\right\}  _{k=1}^{K_{N}}\quad\overset
{(5)}{\Longrightarrow}\quad g_{i,0}(x)\quad\overset{(6)}{\Longrightarrow}\quad
q_{i}(x).
\]

Note that after step (1) the problem is reduced to separate problems on the
edges. The two spectra problem arising after step (2) can be solved by
different existing methods, nevertheless here we propose a method which uses
the fact that the constants $\left\{  g_{i,n}(L_{i})\right\}  _{n=0}^{N}$ and
$\left\{  s_{i,n}(L_{i})\right\}  _{n=0}^{N}$ are also known.

\section{Numerical examples\label{Sect Numerical}}

\textbf{Example 1. }Let us consider a star graph of nine edges of lengths
\begin{equation}
L_{1}=\frac{e}{2},\,L_{2}=1,\,L_{3}=\frac{\pi}{2},\,L_{4}=\frac{\pi}%
{3},\,L_{5}=\frac{e^{2}}{4},\,L_{6}=1.1,\,L_{7}=1.2,\,L_{8}=1,\,L_{9}=1.4
\label{L}%
\end{equation}
The corresponding nine components of the potential are defined as follows%
\[
q_{1}(x)=\left\vert x-1\right\vert +1,\,q_{2}(x)=e^{-(x-\frac{1}{2})^{2}%
},\,q_{3}(x)=\sin\left(  8x\right)  +\frac{2\pi}{3},\,q_{4}(x)=\cos\left(
9x^{2}\right)  +2,
\]%
\[
q_{5}(x)=\frac{1}{x+0.1},\,q_{6}(x)=\frac{1}{\left(  x+0.1\right)  ^{2}%
},\,q_{7}(x)=e^{x},
\]%
\[
\,q_{8}(x)=\left\{
\begin{tabular}
[c]{ll}%
$-35.2x^{2}+17.6x,$ & $0\leq x<0.25$\\
$35.2x^{2}-35.2x+8.8,$ & $0.25\leq x<0.75$\\
$-35.2x^{2}+52.8x-17.6,$ & $0.75\leq x\leq1,$%
\end{tabular}
\ \ \right.  \,q_{9}(x)=J_{0}(9x).
\]

The potential $q_{8}(x)$ (from \cite{Brown et al 2003}, \cite{Rundell Sacks})
is referred to below as saddle potential.

In Fig. 1 we show the exact potentials (continuous line) together with the
recovered ones (marked with asterisks), computed with the proposed algorithm.
Here two elements from each row of the Weyl matrix (that is, $M_{k}=0$)
$\mathbf{M}_{ii}(\rho_{k}^{2})$ and $\mathbf{M}_{ii+1}(\rho_{k}^{2})$ (for
$i=M$, $i+1$ was replaced by $1$) were given at 190 points $\rho_{k}$, and $N$
was chosen as $N=9$ (ten coefficients $g_{n}(L_{i})$ and ten coefficients
$s_{n}(L_{i})$).

\medskip%

%TCIMACRO{\FRAME{ftbpFU}{7.3874in}{3.5197in}{0pt}{\Qcb{The potential of the
%quantum star graph from Example 1 is recovered from two elements from each row
%of the Weyl matrix (that is, $M_{k}=0$) $\QTR{bf}{M}_{ii}(\rho_{k}^{2})$ and
%$\QTR{bf}{M}_{ii+1}(\rho_{k}^{2})$ (for $i=M$, $i+1$ was replaced by $1$)
%given at 190 points $\rho_{k}$ distributed uniformly on the segment
%$[1+0.1i,100+0.1i]$. Here $N=9$. }}{\Qlb{Fig1}}{9edges.eps}%
%width 7.3874in;  height 3.5197in;  depth 0pt;  original-width 13.4509in;
%original-height 6.3714in;  cropleft "0";  croptop "1";  cropright "1";
%cropbottom "0";  filename '9edges.eps';file-properties "XNPEU";}}}%
%BeginExpansion
\begin{figure}
[ptb]
\begin{center}
\includegraphics[
height=3.5197in,
width=7.3874in
]%
{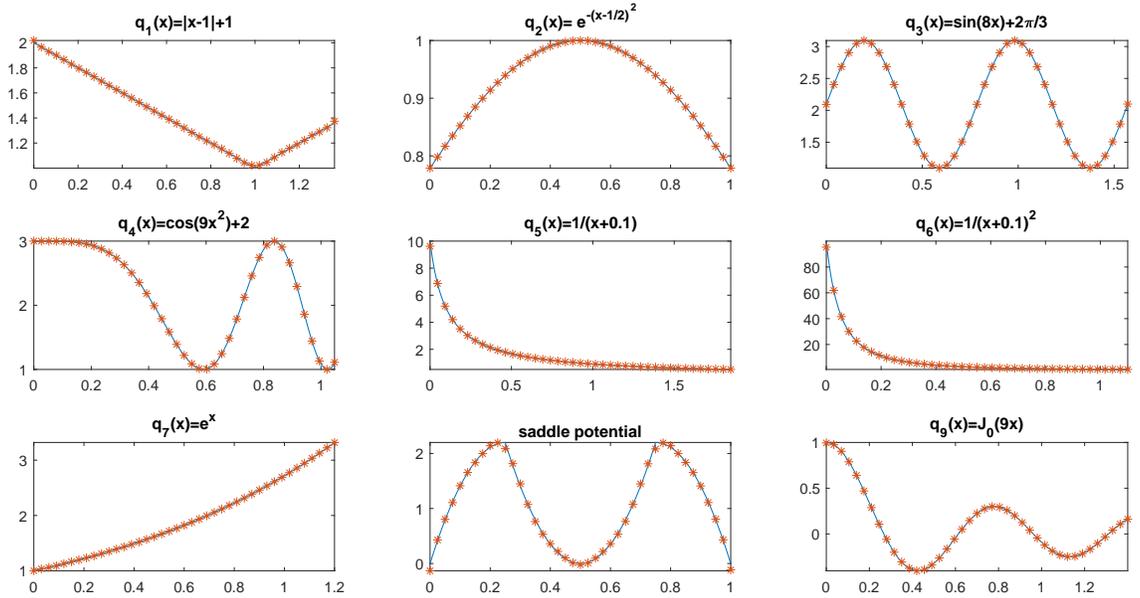}%
\caption{The potential of the quantum star graph from Example 1 is recovered
from two elements from each row of the Weyl matrix (that is, $M_{k}=0$)
$\mathbf{M}_{ii}(\rho_{k}^{2})$ and $\mathbf{M}_{ii+1}(\rho_{k}^{2})$ (for
$i=M$, $i+1$ was replaced by $1$) given at 190 points $\rho_{k}$ distributed
uniformly on the segment $[1+0.1i,100+0.1i]$. Here $N=9$. }%
\label{Fig1}%
\end{center}
\end{figure}
%EndExpansion

Here the maximum relative error was attained in the case of the potential
$q_{6}(x)$, and it resulted in approximately $0.047$ at the endpoint $x=0$.
All other potentials were computed more accurately.

It is interesting to track how accurately the two spectra were computed on
each edge. For example, Table 1 presents some of the \textquotedblleft
exact\textquotedblright\ Dirichlet-Dirichlet eigenvalues on $e_{6}$ computed
with the aid of the Matslise package \cite{Ledoux et al} (first column), the
approximate eigenvalues, computed as described in Subsection
\ref{Subsect reduction to two spectra} by calculating zeros of $S_{6,9}%
(\rho,L_{6})$ and the absolute error of each presented eigenvalue. Notice that
both the absolute and relative errors remain small even for large indices.

\medskip

\bigskip%
\begin{tabular}
[c]{|l|l|l|l|}\hline
\multicolumn{4}{|l|}{Table 1: Dirichlet-Dirichlet eigenvalues of $q_{6}(x)$%
}\\\hline
$n$ & $\lambda_{n}$ & $\widetilde{\lambda}_{n}$ & $\left\vert \lambda
_{n}-\widetilde{\lambda}_{n}\right\vert $\\\hline
$1$ & $11.3620706$ & $11.3620710$ & $4\cdot10^{-7}$\\\hline
$11$ & $994.949643$ & $994.949630$ & $1.3\cdot10^{-5}$\\\hline
$51$ & $21223.885957$ & $21223.885873$ & $8.4\cdot10^{-5}$\\\hline
$101$ & $83214.803376$ & $83214.803222$ & $1.5\cdot10^{-4}$\\\hline
\end{tabular}

\bigskip

Similar results are obtained for the Neumann-Dirichlet spectrum, as shown in
Table 2.

\bigskip%
\begin{tabular}
[c]{|l|l|l|l|}\hline
\multicolumn{4}{|l|}{Table 2: Neumann-Dirichlet eigenvalues of $q_{6}(x)$%
}\\\hline
$n$ & $\lambda_{n}$ & $\widetilde{\lambda}_{n}$ & $\left\vert \lambda
_{n}-\widetilde{\lambda}_{n}\right\vert $\\\hline
$1$ & $10.21124706$ & $10.21124734$ & $2.8\cdot10^{-7}$\\\hline
$11$ & $908.123501$ & $908.123578$ & $7.7\cdot10^{-5}$\\\hline
$51$ & $20809.976547$ & $20809.976583$ & $3.6\cdot10^{-5}$\\\hline
$101$ & $82393.027033$ & $82393.027162$ & $1.3\cdot10^{-4}$\\\hline
\end{tabular}

\bigskip

The number and the distribution of the points $\rho_{k}$ influence the
possibility of an accurate recovery of the potential. So much better results
are obtained when a uniform distribution of the points $\rho_{k}$ is replaced,
e.g., by their logarithmically uniform distribution. For example, the same
accuracy as reported above was obtained for 90 points $\rho_{k}$ chosen
according to the rule $\rho_{k}=10^{\alpha_{k}}+0.1i$ with $\alpha_{k}$ being
uniformly distributed on $\left[  0,2\right]  $. Such choice delivers a set of
points which are more densely distributed near $\rho=1+0.1i$ and more sparsely
near $\rho=100+0.1i$. A shift of the points $\rho_{k}$ to a relatively large
distance from zero leads to a deterioration of the results. For example, when
$\alpha_{k}$ are chosen in $\left[  1,2\right]  $ the method with all the same
parameters fails, while for $\alpha_{k}$ chosen in $\left[  0.5,2\right]  $ it
delivers accurate results. This deterioration of the results may happen due to
the lack of information near the first eigenvalues of the Dirichlet spectrum
of the graph.

Another interesting feature of the method is illustrated by Fig. 2. Namely,
when the number of the points $\rho_{k}$ is small, the use of more equations
of the form (\ref{type2}), i.e., of more elements of the Weyl matrix may help
to improve the accuracy. Here we present the component of the potential
$q_{8}(x)$. The whole potential on the graph was recovered with $N=7$ from 30
points $\rho_{k}$ distributed logarithmically uniformly on the same segment as
above ($\alpha_{k}\in\left[  0,2\right]  $). The best accuracy corresponds to
$M_{k}=7$, that is when all the elements of the Weyl matrix were used, while
the worst result corresponds to $M_{k}=0$, that is when only two elements from
each row were used.%

%TCIMACRO{\FRAME{ftbpFU}{5.6066in}{2.5482in}{0pt}{\Qcb{The component of the
%potential $q_{8}(x)$. The whole potential on the graph was recovered with
%$N=7$ from 30 points $\rho_{k}$ distributed logarithmically uniformly on the
%same segment as above. The best accuracy corresponds to $M_{k}=7$, that is
%when all the elements of the Weyl matrix were used, while the worst result
%corresponds to $M_{k}=0$, that is when only two elements from each row were
%used.}}{}{q8.eps}{\special{ language "Scientific Word";  type "GRAPHIC";
%width 5.6066in;  height 2.5482in;  depth 0pt;  original-width 13.0437in;
%original-height 5.8935in;  cropleft "0";  croptop "1";  cropright "1";
%cropbottom "0";  filename 'q8.eps';file-properties "XNPEU";}}}%
%BeginExpansion
\begin{figure}
[ptb]
\begin{center}
\includegraphics[
height=2.5482in,
width=5.6066in
]%
{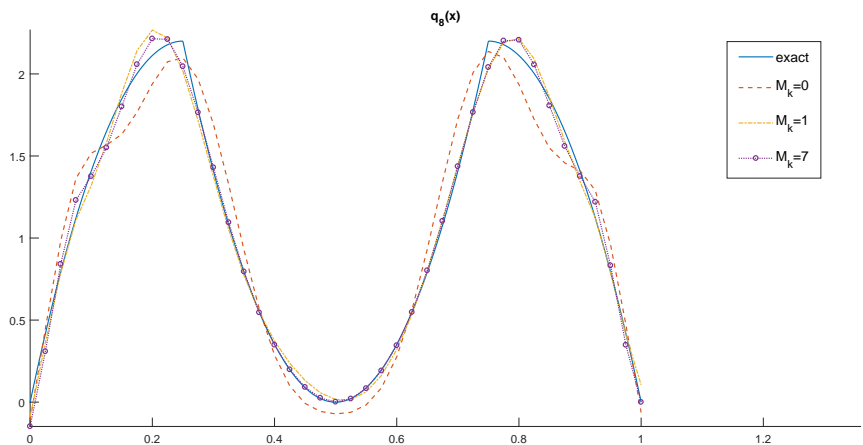}%
\caption{The component of the potential $q_{8}(x)$. The whole potential on the
graph was recovered with $N=7$ from 30 points $\rho_{k}$ distributed
logarithmically uniformly on the same segment as above. The best accuracy
corresponds to $M_{k}=7$, that is when all the elements of the Weyl matrix
were used, while the worst result corresponds to $M_{k}=0$, that is when only
two elements from each row were used.}%
\end{center}
\end{figure}
%EndExpansion

It is worth noting that the whole computation takes few seconds performed in
Matlab 2017 on a Laptop equipped with a Core i7 Intel processor.

The method copes equally well with inverse problems on star graphs with a
larger number of edges, though to obtain a similar accuracy in this case more
input data (points $\rho_{k}$) are required.

\section{Conclusions\label{Sect Concl}}

A new method for solving the inverse problem on quantum star graphs consisting
in the recovery of the potential from the Weyl matrix is developed. The main
role in the proposed approach is played by the coefficients of the Neumann
series of Bessel functions expansion of solutions of the Sturm-Liouville
equation. With their aid the given data lead to separate two-spectra inverse
Sturm-Liouville problems on each edge. These two-spectra problems are solved
by a direct method reducing each problem to a system of linear algebraic
equations, and the crucial observation is that the potential is recovered from
the first component of the solution vector.

The method is simple, direct and accurate. Its performance is illustrated by
numerical examples. In subsequent works we plan to extend this method to other
types of inverse spectral problems and to more general graphs.

\textbf{Funding }The research of Sergei Avdonin was supported in part by the
National Science Foundation, grant DMS 1909869, and by Moscow Center for
Fundamental and Applied Mathematics. The research of Vladislav Kravchenko was
supported by CONACYT, Mexico via the project 284470 and partially performed at
the Regional mathematical center of the Southern Federal University with the
support of the Ministry of Science and Higher Education of Russia, agreement 075-02-2022-893.

\textbf{Data availability} The data that support the findings of this study
are available upon reasonable request.

\textbf{Declarations}

\textbf{Conflict of interest} The authors declare no competing interests.

\end{document}